\documentclass[12pt]{article} 
\usepackage[sectionbib]{natbib}
\usepackage{array,epsfig,fancyheadings,rotating}
\usepackage[]{hyperref}  
\usepackage[usenames,dvipsnames,svgnames,table]{xcolor} 
\usepackage{sectsty, secdot}
\sectionfont{\fontsize{12}{14pt plus.8pt minus .6pt}\selectfont}
\renewcommand{\theequation}{\thesection\arabic{equation}}
\subsectionfont{\fontsize{12}{14pt plus.8pt minus .6pt}\selectfont}

\usepackage{amsmath}
\usepackage{amssymb}
\usepackage{amsfonts}
\usepackage{multirow}
\usepackage{amsthm}
\usepackage{mathrsfs}
\usepackage{cleveref}

\newtheorem{theorem}{Theorem}
\newtheorem{lemma}{Lemma}
\newtheorem{corollary}{Corollary}
\newtheorem{Assumption}{Assumption}

\theoremstyle{definition}
\newtheorem{definition}{Definition}
\newtheorem{condition}{Condition}

\newtheorem{remark}{Remark}

\newcommand{\bea}{\begin{eqnarray*}}
	\newcommand{\eea}{\end{eqnarray*}}
\newcommand{\bean}{\begin{eqnarray}}
\newcommand{\eean}{\end{eqnarray}}
\newcommand{\baln}{\begin{align}}
\newcommand{\ealn}{\end{align}}
\newcommand{\bdefi}{\begin{definition}}
	\newcommand{\edefi}{\end{definition}}
\newcommand{\bi}{\begin{itemize}}
	\newcommand{\ei}{\end{itemize}}
\newcommand{\benu}{\begin{enumerate}}
	\newcommand{\eenu}{\end{enumerate}}
\newcommand{\ben}{\begin{enumerate}}
	\newcommand{\een}{\end{enumerate}}

\DeclareMathOperator*{\argmax}{argmax}

\newcommand{\sg}{\Sigma}

\newcommand{\bbR}{\mathbb{R}}

\begin{document}


\newcommand{\bch}{\color{blue}  }   
\newcommand{\ech}{\color{black}  ~}    

\renewcommand{\baselinestretch}{2}

\markright{ \hbox{\footnotesize\rm Statistica Sinica
}\hfill\\[-13pt]
\hbox{\footnotesize\rm
}\hfill }

\markboth{\hfill{\footnotesize\rm XUAN CAO AND KYOUNGJAE LEE} \hfill}
{\hfill {\footnotesize\rm BAYESIAN JOINT SELECTION} \hfill}

\renewcommand{\thefootnote}{} 
$\ $\par


\fontsize{12}{14pt plus.8pt minus .6pt}\selectfont \vspace{0.8pc}

\centerline{\large\bf JOINT BAYESIAN VARIABLE AND DAG SELECTION CONSISTENCY  }
\vspace{2pt} \centerline{\large\bf FOR HIGH-DIMENSIONAL REGRESSION MODELS}
\vspace{2pt} \centerline{\large\bf WITH NETWORK-STRUCTURED COVARIATES}
\vspace{.4cm} \centerline{Xuan Cao and Kyoungjae Lee} \vspace{.4cm} \centerline{\it
	University of Cincinnati and Inha University} \vspace{.55cm} \fontsize{9}{11.5pt plus.8pt minus
	.6pt}\selectfont


\begin{quotation}
\noindent {\it Abstract:}
We consider the joint sparse estimation of regression coefficients and the covariance matrix for covariates in a high-dimensional regression model, where the predictors are both relevant to a response variable of interest and  functionally related to one another via a Gaussian directed acyclic graph (DAG) model. Gaussian DAG models introduce sparsity in the Cholesky factor of the inverse covariance matrix, and the sparsity pattern in turn corresponds to specific conditional independence assumptions on the underlying predictors. A variety of methods have been developed in recent years for Bayesian inference in identifying such network-structured predictors in regression setting, yet crucial sparsity selection properties for these models have not been thoroughly investigated. 
In this paper, we consider a hierarchical model with spike and slab priors on the regression coefficients and a flexible and general class of DAG-Wishart distributions with multiple shape parameters on the Cholesky factors of the inverse covariance matrix. Under mild regularity assumptions, we establish the joint selection consistency for both the variable and the underlying DAG of the covariates when the dimension of predictors is allowed to grow much larger than the sample size. We demonstrate that our method outperforms existing methods in selecting network-structured predictors in several simulation settings.

\vspace{9pt}
\par
\end{quotation}\par

\def\thefigure{\arabic{figure}}
\def\thetable{\arabic{table}}

\renewcommand{\theequation}{\thesection.\arabic{equation}}

\fontsize{12}{14pt plus.8pt minus .6pt}\selectfont

\setcounter{equation}{0} 
\section{Introduction}
In modern day statistics, datasets where the number of variables is much larger than the number of samples are more pervasive than they have ever been. 
One of the major problems is high-dimensional variable selection, where the challenge is to select a subset of predictor variables which significantly affect a given response. The literature on Bayesian variable selection in linear regression is vast and rich. 
\citet{George:McCulloch:1993} propose the stochastic search variable selection which uses the Gaussian distribution with a zero mean and a small but fixed variance as the spike prior, and another Gaussian distribution with a large variance as the slab prior. \citet*{ishwaran2005} also use Gaussian spike and slab priors, but with continuous bimodal priors for the variance of the regression coefficient to alleviate the difficulty of choosing specific prior parameters. \citet*{Narisetty:He:2014} introduce shrinking and diffusing priors as spike and slab priors, and establish model selection consistency of the approach in a high-dimensional setting. 

Another important problem is how to formulate models and develop inferential procedures to understand the complex relationships and multivariate dependencies in these high-dimensional datasets. 
A covariance matrix is one of the most fundamental objects that quantifies these relationships. A common and effective approach for covariance estimation in sample starved settings is to induce sparsity either in the covariance matrix, its inverse, or the Cholesky factor of the inverse. The sparsity patterns in these matrices can be uniquely encoded in terms of appropriate graphs. Hence the corresponding models are often referred to as covariance graph models (sparsity in $\Sigma$), concentration graph models (sparsity in $\Omega = \Sigma^{-1}$), and directed acyclic graph (DAG) models (sparsity in the Cholesky factor of $\Omega$). 

In this work, we will work in a high-dimensional regression setting, where the predictors are both relevant to a response variable of interest and  functionally related to one another via a Gaussian DAG model.
Our goal is to jointly perform variable and DAG selection, and to establish the selection consistency in a high-dimensional regime. 
The advantage of joint modeling is that we can borrow information from the DAG structure to improve the performance of variable selection.
One popular motivation for this type of problem comes from genomic studies: the mechanism for an effect on an outcome such as a quantitative molecular phenotypes including gene expression, proteomics, or metabolomics data often displays a coordinated change along a pathway, and the impact of one single genotype may not be apparent. 
In this setting, our proposed method can incorporate and highlight unknown pathways or regulatory networks that impact the response, which can potentially improve the performance of variable selection by borrowing information from the network structure. 
To uncover these relationships, we develop a Bayesian hierarchical model that favors the inclusion of variables that are not only relevant to the outcome of interest but also linked through a DAG. 

When the underlying graph structure is known, several approaches including both frequentist and Bayesian methods have been proposed and studied in the literature to solve the variable selection problem. \citet{li2008network, Li:Li:2010} study a graph-constrained regularization procedure and its theoretical properties to take into account the neighborhood information of the variables measured on a known graph. \citet{Pan:Xie:Shen:2010} propose a grouped penalty based on the $L_\gamma$-norm that smooths the regression coefficients of the predictors over the available network. On the Bayesian side, \citet{Li:Zhang:2010} and \citet{Stingo:Vannucci:2010} incorporate a graph structure in the Markov random field (MRF) prior on indicators of variable selection, encouraging the joint selection of predictors with known relationships. \citet{Stingo:Chen:Vannucci:2011} and \citet{Peng:2013} propose the selection of both pathways and genes within them based on prior knowledge on gene-gene interactions or functional relationships.

However, when the underlying graph is unknown and needs to be selected, comparatively fewer methods have been proposed.  \citet{Dobra:2009} estimate a network among relevant predictors by first performing a stochastic search in the regression setting to identify possible subsets of predictors, then applying a Bayesian model averaging method to estimate a dependency network. \citet{Liu:Chakraborty:2014}  develop a Bayesian method for regularized regression, which provides inference on the inter-relationship between variables by explicitly modeling through a graph Laplacian matrix. \citet{Joint:Network:2016} simultaneously infer a sparse network among the predictors and perform variable selection using this network as guidance by incorporating it into a prior favoring selection of connected variables based on a Gaussian graphical model among the predictors, which provides a sparse and interpretable representation of the conditional dependencies found in the data. 
In a slightly different context, \cite{chekouo2015mirna} and \cite{chekouo2016bayesian} relate two sets of covariates via a DAG to integrate multiple genomic platforms and select the most relevant features. Given the ordering of variables, they use a mixture of a non-local prior \citep{johnson2012bayesian} and a point mass at zero to infer the DAG structure.

Despite the developments in Bayesian methods for joint variable and graph selection, a rigorous investigation of the high-dimensional consistency properties of these methods has not been undertaken to the best of our knowledge. Hence, our goal was to investigate if joint selection consistency results could be established in the high-dimensional  regression setting with network-structured predictors. This is a challenging goal, particularly because of the interaction between the regression coefficients and the graph in the posterior analysis, and the massive parameter space to be explored for both the coefficients and the graph.  

In this paper, we consider a hierarchical multivariate regression model with DAG-Wishart priors on the covariance matrix for the predictors, spike and slab priors on regression coefficients, independent Bernoulli priors for each edge in the DAG, and a MRF prior linking the variable indicators to the graph structure. 
Under high-dimensional settings, we establish {posterior ratio consistency}, following the nomenclature in \citet{CKG:2017} and \citet{Narisetty:He:2014}, for both the variable and the DAG with given DAG and variable, respectively (Theorems \ref{thm1} and \ref{thm2}). 
In Theorems \ref{thm3} and \ref{thm4},  we also establish the posterior ratio consistency and the strong selection consistency for any pair of the DAG and variable. In particular, the strong selection consistency implies that under the true model, the posterior probability of the true variable indicator and the true graph converges in probability to $1$ as $n \rightarrow \infty$. Finally, through simulation studies, we demonstrate that the  models studied in this paper can outperform existing state-of-the-art methods including both penalized likelihood and Bayesian approaches in several settings.

The rest of paper is organized as follows. Section \ref{sec2} provides background material regarding Gaussian DAG model and the DAG-Wishart distribution. In Section \ref{sec:modelspecification}, we introduce our hierarchical Bayesian model. Model selection consistency results are stated in Section \ref{sec:modelconsistency} with proofs provided in Supplementary material. 
In Section \ref{sec:experiments}, we conduct simulation experiments to illustrate the performance of the proposed method. 
Benefits of our Bayesian method for identifying  network-structured predictors are demonstrated vis-a-vis existing Bayesian and penalized likelihood approaches. 
We end our paper with a discussion session in Section \ref{sec:discussion}.

\setcounter{equation}{0} 
\section{Preliminaries}\label{sec2}

\noindent
In this section, we provide the necessary background material from graph theory, 
Gaussian DAG models, and DAG-Wishart distributions. 

\subsection{Gaussian DAG models} \label{sec2.1}

\noindent
Throughout this paper, a directed acyclic graph (DAG) $\mathscr{D} = (V,E)$ consists 
of a vertex set $V = \{1,\ldots,p\}$ and an edge set $E$ such that there is no directed 
path starting and ending at the same vertex. As in \cite{BLMR:2016} and \cite{CKG:2017}, we will assume a parent ordering, where that all the edges are directed from larger vertices to smaller vertices. 
Thus, the ordering of variables is assumed to be known throughout the paper. 
The set of parents of $i$, denoted by $pa_i(\mathscr D)$, is the collection of all vertices which are larger than $i$ and share an edge with $i$. A Gaussian DAG model over a given DAG $\mathscr{D}$, denoted by 
$\mathscr{N}_{\mathscr{D}}$, consists of all multivariate Gaussian distributions which 
obey the directed Markov property with respect to a DAG $\mathscr{D}$. In particular, 
if $ {x}=(x_1, \ldots, x_p)^T \sim N_p(0,\Sigma)$ and $N_p(0,\Sigma) \in 
\mathscr{N}_{\mathscr{D}}$, then $x_i \perp  {x}_{\{i+1,\ldots,p\}\backslash pa_i(\mathscr D)}|
{x}_{pa_i(\mathscr D)}$ for each $i$. 

Any positive definite matrix $\Omega$ can be uniquely decomposed as $\Omega = 
LD^{-1}L^T$, where $L$ is a lower triangular matrix with unit diagonal entries, and $D$ 
is a diagonal matrix with positive diagonal entries. This decomposition is known as the 
modified Cholesky decomposition of $\Omega$ (see for example
\citet{Pourahmadi:2007}). It is well-known that if $\Omega = LD^{-1}L^T$ is the 
modified Cholesky decomposition of $\Omega$, then $N_p(0,\Omega^{-1}) \in 
\mathscr{N}_{\mathscr{D}}$ if and only if $L_{ij} = 0$ whenever $i \notin pa_j (\mathscr D)$. In other 
words, the structure of the DAG $\mathscr{D}$ is reflected in the Cholesky factor $L$ of 
the inverse covariance matrix. 

Given a DAG $\mathscr{D}$ on $p$ vertices, denote $\mathscr{L}_{\mathscr{D}}$ as 
the set of lower triangular matrices with unit diagonals and $L_{ij} = 0$ if $i \notin 
pa_j(\mathscr D)$, and let $\mathscr{D}_+^p$ be the set of strictly positive diagonal matrices in 
$\mathbb{R}^{p \times p}$. We refer to $\Theta_{\mathscr{D}} = \mathscr{D}_+^p 
\times \mathscr{L}_{\mathscr{D}}$ as the Cholesky space corresponding to 
$\mathscr{D}$, and $(D,L) \in \Theta_{\mathscr{D}}$ as the Cholesky parameter 
corresponding to $\mathscr{D}$. In fact, the relationship between the DAG and the 
Cholesky parameter implies that 
$
\mathscr{N}_{\mathscr{D}} = \{N_p(0,(L^T)^{-1}DL^{-1}):(D,L) \in 
\Theta_{\mathscr{D}}\}. 
$

The skeleton of $\mathscr{D}$, denoted by $\mathscr{D}^u = (V,E^u)$, can be 
obtained by replacing all the directed edges of $\mathscr{D}$ by undirected ones. We define the adjacency matrix of $\mathscr D$ to be a (0,1)-matrix such that the elements of the matrix indicate whether pairs of vertices are adjacent or not in $\mathscr D$, i.e., 1 representing adjacent, 0 representing not adjacent.

\subsection{DAG-Wishart Distribution} \label{sec2.2}

\noindent
In this section, we revisit the multiple shape parameter DAG-Wishart distributions 
introduced in \citet{BLMR:2016}. 
Given a directed 
graph $\mathscr{D} = (V,E)$ with $V = \{1, \ldots, p\}$ and a $p \times p$ matrix $A$, 
denote the column vectors $A_{\mathscr D .i}^> = (A_{ij})^T_{j \in pa_{i}(\mathscr D)}$ and 
$A_{\mathscr D.i}^{\ge} = (A_{ii}, (A_{\mathscr D.i}^>)^T)^T.$ 
Also,  
$$ A_{\mathscr D}^{ \ge i} = 
\left[ \begin{matrix}
A_{ii} & (A_{\mathscr D.i}^>)^T \\
A_{\mathscr D.i}^> & A_{\mathscr D}^{>i}
\end{matrix} \right] ,
$$
where $A_{\mathscr D}^{>i} = (A_{kj})_{k,j \in pa_{i}(\mathscr D)}$.
In particular, we have $A_{\mathscr D.p}^{\ge} = A_{\mathscr D}^{ \ge p} = A_{pp}$. 
Let $\nu_{i}(\mathscr D) = |pa_{i} (\mathscr{D})| = |\{j: j>i, (j,i) \in E({\mathscr{D}})\}|$. 

The DAG-Wishart distributions in \citet{BLMR:2016} corresponding to a DAG 
$\mathscr{D}$ are defined on the Cholesky space $\Theta_{\mathscr{D}}$. 
Given a $p \times p$ positive definite matrix $U$ and a $p$-dimensional vector ${\boldsymbol \alpha} (\mathscr{D}) = ( \alpha_1 (\mathscr{D}), \ldots, \alpha_p (\mathscr{D}) )$ with $\min_{1\le i \le p} \{\alpha_i(\mathscr D) - \nu_i(\mathscr D)\} >2$, the probability density of the DAG-Wishart distribution is given by 
\begin{equation} \label{a2}
\pi_{U,  \alpha(\mathscr D)}^{\Theta_{\mathscr{D}}} (D,L) = 
\frac{1}{z_{\mathscr{D}}(U,{\boldsymbol \alpha} (\mathscr{D}))} \exp\{-\frac12\mbox{tr}
((LD^{-1}L^T)U)\} \prod_{i=1}^p D_{ii}^{-\frac{\alpha_i (\mathscr{D})}2} \mathbb I \Big( (D,L) \in \Theta_{\mathscr{D}}\Big) ,
\end{equation}
where 
$$z_{\mathscr{D}}(U,{\boldsymbol \alpha} (\mathscr{D})) = \prod_{i=1}^{p} \frac{\Gamma(\frac{\alpha_i (\mathscr{D})}2 - \frac{\nu_{i}(\mathscr D)}2 - 	1)2^{\frac{\alpha_i (\mathscr{D})}2 - 1}(\sqrt{\pi})^{\nu_{i}(\mathscr D)} 	det(U_{\mathscr D}^{>i})^{\frac{\alpha_i (\mathscr{D})}2 - \frac{\nu_{i}(\mathscr D)}2 - 		\frac32}}{det(U_{\mathscr D}^{\ge i})^{\frac{\alpha_i (\mathscr{D})}2 - \frac{\nu_{i} (\mathscr D)}2 - 1}}$$
and $\mathbb{I}(\cdot)$ stands for the indicator function.
The above density has the same form as the classical Wishart density, but is defined on the lower dimensional space  $\Theta_{\mathscr{D}}$ and has $p$ shape parameters $\{\alpha_i (\mathscr{D})\}_{i=1}^p$ which can be used for differential shrinkage of  variables in high-dimensional settings. 

The class of densities $\pi_{U,  \alpha(\mathscr D)}^{\Theta_{\mathscr{D}}}$ form a 
conjugate family of priors for the Gaussian DAG model $\mathscr{N}(\mathscr{D})$. 
In particular, if the prior on $(D,L) \in \Theta_{\mathscr{D}}$ is $\pi_{U,  \alpha(\mathscr D)}^{\Theta_{\mathscr{D}}}$ and $ {X}_1, \ldots,  {X}_n \mid D,L, \mathscr{D} \overset{i.i.d.}{\sim} N_p({0},(L^T)^{-1}DL^{-1})$, then the resulting posterior distribution of $(D,L)$ is $\pi_{\tilde{U},\tilde{  \alpha}(\mathscr D)}^{\Theta_{\mathscr{D}}}$, where $S = \frac1n\sum_{i=1}^n {X}_i {X}_i^T$, $\tilde{U} = U + nS$, and $\tilde{\alpha}
(\mathscr D) = (n+\alpha_1(\mathscr D), \ldots, n+\alpha_p(\mathscr D))$.

\setcounter{equation}{0} 
\section{Model Specification} \label{sec:modelspecification}
In this section, we specify our hierarchical model to facilitate the purpose of joint variable and DAG selection for regression models with network-structured predictors. We start by considering the standard Gaussian linear regression model with $p$ coefficients and by introducing some required 
notations. Similar to \citet{Joint:Network:2016} and \citet{li2008network}, consider both the response $  Y = (y_1, \ldots, y_n)\in \mathbb{R}^{n\times 1}$ and the predictors $X = \left(  X_1,    \ldots,   X_n\right)^T \in \mathbb{R}^{n\times p}$ to be random variables. In particular, $  Y \sim N_n \left(X \beta, \sigma^2 I_n\right)$, and the predictors are assumed to obey a multivariate Gaussian distribution, i.e., $  X_i \overset{i.i.d.}{\sim} N_p \left(0, (LD^{-1}L^T)^{-1}\right)$ for $i = 1,2, \ldots, n$, where $\beta \in  \bbR^{p \times 1}$ is a vector of regression coefficients and $(L,D)$ represents the Cholesky parameter corresponding to a DAG $\mathscr D$. Let symmetric matrix $G=\left(G_{ij}\right)_{1\le i,j\le p}$ represent the adjacency matrix corresponding to DAG $\mathscr D$ where $G_{ij} = G_{ji} = 1$ if and only if there is an edge between vertex $i$ and vertex $j$, and $G_{ij} = G_{ji} = 0$ otherwise. Our goal is both (i) the variable selection, i.e., to correctly identify all the non-zero regression  
coefficients, and  (ii) network estimation, i.e., to precisely recover the sparsity pattern in $\mathscr D$. 

For variable selection, we denote a variable indicator $\gamma = \left\{\gamma_1,  \ldots, \gamma_p\right\}$,  where $\gamma_j\in \{0,1\}$ for $1\le j \le p$, and $|\gamma|= \sum_{j=1}^p \gamma_j$. Let $  \beta_\gamma = ( \beta_j )^T_{\{j: \gamma_j =1\}}\in \bbR^{|\gamma|\times 1}$  be the vector formed by the active components in $\beta$ corresponding to a model $\gamma$. 
For any $n \times p$ matrix $A$, let $A_k$ represent the submatrix formed from the columns of $A$ corresponding to model $k$. 
In particular, Let $X_\gamma$ denote the design matrix formed from the columns of 
$  X$ corresponding to model $\gamma$. 
For the network estimation, the class of DAG-Wishart distribution in Section \ref{sec2.2} can be used for joint variable and DAG selection through the following hierarchical model.
\begin{eqnarray} 
&      Y |   X_\gamma, \beta_\gamma \sim N_n \left(X_\gamma \beta_\gamma, \sigma^2 I_n\right), \label{model:spec1}
\\
&    X_i | (L, D), \mathscr D \overset{i.i.d.}{\sim} N_p \left(0, (LD^{-1}L^T)^{-1}\right), \quad \mbox{for } i = 1, 2, \ldots, n, \label{model:spec2}
\\
&(L, D) | \mathscr D \sim \pi_{U,  \alpha(\mathscr D)}^{\Theta_{\mathscr{D}}} (D,L), \label{model:spec3}
\\
& \beta_\gamma |   \gamma \sim N_{|\gamma|} \left(0, \tau^2\sigma^2 I_\gamma\right), \label{model:spec4}
\\
& \pi(\mathscr D) \propto \prod_{j = 1}^{p-1} q^{\nu_j(\mathscr D)}(1-q)^{p-j-\nu_j(\mathscr D)}\mathbb I\left\{\max_{1 \le j \le p-1}\nu_i(\mathscr D) < R \right\}, \label{model:spec5}
\\
& \pi(  \gamma | \mathscr D) \propto \exp\left(-a 1^T \gamma + b \gamma^T G \gamma \right)\mathbb I\left\{|\gamma| <  R\right\}. \label{model:spec6}
\end{eqnarray}  
for some  constants $\sigma, \tau, a >0$, $b \ge 0$, $0<q<1$ and a positive integer $0\le R \le p$. Here we assume that $\sigma$ in \eqref{model:spec1} is a known constant for simplicity.
However, it can be  extended to unknown $\sigma$ case by imposing an inverse-gamma prior, which will be shown in Corollary \ref{cor1}.
Note that in (\ref{model:spec4}), we are essentially imposing a spike and slab prior on the regression coefficients, where $\tau^2$ indicates the variance of the slab part. See \citet{Narisetty:He:2014}, \citet{yang:2016} and the references therein. 
Prior (\ref{model:spec5}) corresponds to an Erdos-Renyi type of prior over the space of DAGs. In particular, similar to \citet{CKG:2017}, define $e_{ji} = \mathbb{I} \{(j,i) \in E(\mathscr D)\}$, $1 \le j < i \le p$ to be the edge indicator. 
Let $e_{ji}$, $1 \le i < j < p$ be independent identically distributed Bernoulli($q$) random variables. 
Recall $\nu_j(\mathscr D) = |pa_j(\mathscr D)|$ is the cardinality of the parent set of vertex $j$. It follows that $\pi (\mathscr D) = \prod_{(j,i):1 \leq j < i \leq p} q^{e_{ji}} \left(1-q\right)^{1-e_{ji}} = \prod_{j=1}^{p-1} q^{\nu_j(\mathscr{D})} (1-q)^{p-j-\nu_j(\mathscr{D})}.$ 
In \eqref{model:spec5} and \eqref{model:spec6}, the positive integer $R$ is an upper bound on the DAG and regression complexity.
Note that to obtain our desired asymptotic consistency results, appropriate conditions for these hyperparameters $\tau, R, a, b$ as well as the edge probability $q$ will be introduced in Section \ref{sec:modelconsistency}. 
\begin{remark}
	In (\ref{model:spec6}), given a DAG $\mathscr{D}$, we are imposing a Markov random field (MRF) prior on the variable indicator $\gamma$ that favors the inclusion of variables linked to other variables in the associated DAG. 
	MRF priors have also been used in the variable selection setting in \citet{Joint:Network:2016, Li:Zhang:2010} and \citet{Stingo:Vannucci:2010}. In particular, as indicated in \citet{Joint:Network:2016}, the parameter $a$ in (\ref{model:spec6}) controls the variable inclusion probability, with  larger values of $a$ corresponding to sparser models, while $b$ essentially determines how strongly the inclusion probability of a variable is affected by the inclusion of its neighbors in the DAG.
\end{remark}

The hierarchical model in (\ref{model:spec1})-(\ref{model:spec6}) can be used to estimate a  pair of variable and DAG as follows. By 
(\ref{a2}) and Bayes' rule, the following lemma gives the (marginal) joint posterior probabilities with proof provided in the Supplementary material.
\begin{lemma} \label{lemma_conditional}
	Under the hierarchical model in (\ref{model:spec1})-(\ref{model:spec6}), the (marginal) joint variable and DAG posterior is given by,
	\begin{align} \label{conditional}
	&\pi\left(\gamma, \mathscr D | Y, X\right) \nonumber \\
	\propto &\,\, \pi(\gamma|\mathscr D)\pi(\mathscr D)\frac{z_{\mathscr{D}}(U+X^TX,n+\alpha(\mathscr D))}{z_{\mathscr{D}}(U, \alpha(\mathscr D))} \nonumber\\
	&\times \det\left(\tau^2X_\gamma^T X_\gamma + I_{|\gamma|} \right)^{-\frac 1 2 }\exp\left\{-\frac{1}{2\sigma^2}\left(Y^T\left(I_n + \tau^2 X_{\gamma} X_{\gamma}^T\right)^{-1}Y\right)\right\},
	\end{align}
	where $z_{\mathscr D}(\cdot, \cdot)$ is the normalized constant in the DAG-Wishart distribution. 	
\end{lemma}
\noindent
Hence, after integrating out $\beta_\gamma$, we have the joint posterior available in closed form (up to the multiplicative constant $\pi(X, Y)$). In particular, these posterior 
probabilities can be used to select a pair of variable and DAG by computing the posterior mode defined by
\begin{equation} \label{a4}
(\hat{\gamma}, \hat{\mathscr D}) =  \argmax_{(\gamma, \mathscr{D})} \pi\left(\gamma, \mathscr D | Y, X\right).
\end{equation}

\setcounter{equation}{0} 
\section{Joint Selection Consistency} \label{sec:modelconsistency}
In this section we will explore the high-dimensional asymptotic properties of the Bayesian joint variable and DAG selection approach specified in Section \ref{sec:modelspecification}. For this purpose, we will work in a setting where the number of regression coefficients $p = p_n$ 
increases with the sample size $n$. The true data generating mechanism is given by 
\begin{eqnarray*}
	Y &=& X \beta_0^n + {\epsilon}_n,
\end{eqnarray*}	
where $Y=(Y_1,\ldots, Y_n)\in \mathbb{R}^n$, $X = (X_1,  \ldots, X_n)^T \in \mathbb{R}^{n\times p_n}$, $X_i \overset{i.i.d.}{\sim} N_{p_n}\left(0,  \Sigma_0^n\right)$ and $\epsilon_n \sim N_n(0, \sigma_0^2 I_n)$.
Here $\beta_0^n$ is the true $p_n$-dimensional vector of regression coefficients, and  $\sg_0^n$ is the true covariance matrix.
As in the usual context of variable selection, we assume that the true vector of regression coefficients is sparse, i.e., all the entries of $\beta_0^n$ are zero except those corresponding to the active entries in the true variable indicator $\gamma_0^n$ \citep{Castillo:2015, yang:2016, Narisetty:He:2014}. 
Denote  ${\rho_1}_n = \min_{j \in \gamma_0^n}|{\beta_0^n}_j|$ and  ${\rho_2}_n = \max_{j \in \gamma_0^n}|{\beta_0^n}_j|$ as the minimum and maximum magnitude of non-zero entries in $\beta_0^n$, respectively. 
We assume that the true quantities $|\gamma_0^n|$, ${\rho_1}_n$ and ${\rho_2}_n$ vary with $n$.
Let $\Omega_0^n = (\Sigma_0^n)^{-1} = L_0^n(D_0^n)^{-1}(L_0^n)^T$, where $(D_0^n, L_0^n)$ denotes the modified Cholesky parameter of 
$\Omega_0^n$. Let $\mathscr D_0^n$ be the true underlying DAG with structure corresponding to the sparsity pattern in $L_0^n$, i.e, $L_0^n \in \mathscr{L}_{\mathscr{D}_0^n}$, and let $G_0^n$ be the adjacency matrix for $\mathscr D_0^n$. 
Denote $d_n$ as the maximum number of non-zero entries in any column of $L_0^n$, and $s_n = \min_{1 \leq j \leq p_n, i \in pa_j(\mathscr{D}_0^n)} |(L_0^n)_{ij}|$  as the minimum magnitude of non-zero off-diagonal entry in $L_0^n$.
Let $\bar P$ denote the probability measure corresponding to the true model presented above. In order to establish the desirable consistency results, we need the following mild assumptions. 
Each assumption is followed by an interpretation/discussion.

\begin{Assumption} \label{assumption:precision}
	There exists $0<\epsilon_{0} \le 1$ such that   $ \epsilon_{0} \le eig_1({\Omega}_0^n) \le eig_{p_n}({\Omega}_0^n) \le \epsilon_{0}^{-1}$ for every $n \ge 1$, where $eig_1({\Omega}_0^n)$ and $eig_{p_n}({\Omega}_0^n)$ are the minimum and maximum eigenvalues of $\Omega_0^n$, respectively.
\end{Assumption}
\noindent
This is a standard assumption for high dimensional covariance asymptotic consistency, both in the frequentist and Bayesian paradigms. See for example \citet{Bickel:Levina:2008, ElKaroui:2008, Banerjee:Ghosal:2014, XKG:2015} and \citet{Banerjee:Ghosal:2015}. \citet{CKG:2017} relax this assumption by allowing the lower and upper bounds on the eigenvalues to depend on $p_n$ and $n$.
\begin{Assumption}  \label{assumption:subexponetial}
	For the true DAG, $d_n\sqrt{ \log p_n /n}\rightarrow 0$ and $d_n\log p_n /(s_n^2n) \rightarrow 0$. For the true regression coefficient,  $|\gamma_0^n|\sqrt{ \log p_n /n} \rightarrow 0$, $\log n\log p_n/(n{\rho_1}_n^2) \rightarrow 0$ and ${\rho_2}_n/\sqrt{\log p_n} \rightarrow 0$ as $n \rightarrow \infty$.
\end{Assumption}
\noindent
This assumption resembles the dimension assumption in \citet{CKG:multi}, and is a much weaker assumption for high dimensional covariance asymptotic than for example, \citet{XKG:2015, Banerjee:Ghosal:2014, Banerjee:Ghosal:2015} and \citet{CKG:2017}. Here we essentially allow the dimension of our covariance matrix to grow slower than $\exp(n/d_n^2)$. 	
Recall that $s_n$ is the smallest (in absolute value) non-zero off-diagonal entry in $L_0^n$, so the second condition in Assumption \ref{assumption:subexponetial} can also be interpreted as the lower bound for the signal size. 
This assumption also known as the ``beta-min" condition provides a lower bound for the signal size that is needed for establishing consistency. 
This type of condition has been used for the exact support recovery of the high-dimensional linear regression models as well as Gaussian DAG models. See for example \citet{yang:2016, KORR:2017,LLL:2018} and \citet{CKG:2017}. 
Assumption \ref{assumption:subexponetial} also allows the complexity of $\gamma_0^n$ as well as the non-zero entries of $\beta_0^n$  to grow with $n$ while stay uniformly bounded by a function of $n$ and $p_n$. In addition, the assumption on ${\rho_1}_n$ can be viewed as the beta-min condition in the regression context.

\begin{Assumption}\label{assumption:hyper:markov}
	The hyperparameters in model (\ref{model:spec4}) and the MRF prior (\ref{model:spec6}) satisfy $\tau^2 \sim \sqrt{\log p_n}$, $a \sim \alpha_1\log p_n$, and $ b n^2/ \{(\log n)^2\log p_n\} \rightarrow 0$ as $n \rightarrow \infty$, where for any positive sequences $a_n$ and $b_n$, $a_n \sim b_n$ implies that there exist positive constants $c$ and $C$ such that $c \le \min(a_n/b_n, b_n/ a_n ) \le \max (a_n/b_n, b_n/ a_n ) \le C$.
\end{Assumption}
\noindent
Recall that the parameter $a$ in (\ref{model:spec6}) controls the variable inclusion probability, and $b$ reflects that how strongly the inclusion probability of a variable is affected by the inclusion of its neighbors in the DAG. In Section \ref{sec:b=0}, we investigate the behavior of the posterior probability evaluated at the true model under $b > 0$ and $b = 0$. In the Bayesian variable selection literature, similar priors corresponding to $a = C \log p_n$ for some constant $C>0$ and $b=0$ have been commonly used to obtain selection consistency \citep{Narisetty:He:2014, Castillo:2015, yang:2016}. The assumption on the the variance of the slab prior, $\tau^2$, is required to approach infinity is also stated here to ensure desired model selection consistency.

\begin{Assumption} \label{assumption:q}
	Let $q_n = O(p_n^{-\alpha_1})$ for some constant $\alpha_1 > 0$ and $R_n$ in model (\ref{model:spec5}) and (\ref{model:spec6}) satisfy $R_n \sim n/\log n$ and $bR_n^2/\log p_n \rightarrow 0$ as $n \rightarrow \infty$.
\end{Assumption}
\noindent
This assumption provides the rate at which the edge probability $q_n$ needs to approach zero. 
It also states that the prior on the space of the $2^{\binom{p_n}{2}}$ possible models, places zero mass on unrealistically large models. Note that  $q_n$ is of slower rate approaching zero compared to the one in \citet{CKG:2017}, which helps avoiding the potential computation limitation such as simulation results always favor the most sparse model. 
This assumption also states that the MRF prior on the space of the $2^{p_n}$ possible models, places zero mass on unrealistically large models (see similar assumptions in \citet{Shin.M:2015,Narisetty:He:2014} in the context of regression).
\begin{Assumption} \label{assumption:hyper}
	For every $n \ge 1$, the hyperparameters for the DAG-Wishart prior 
	$\pi_{U_n,  {\alpha}(\mathscr{D}_n)}^{\Theta_{\mathscr{D}_n}}$ satisfy 
	(i) $2 < \alpha_i(\mathscr{D}_n) - \nu_i(\mathscr D_n) < c$ for every $\mathscr{D}_n$ 
	and $1 \le i \le q_n$, and (ii) $0 < \delta_1 \le eig_1(U_n) \le eig_{p}(U_n) \le 
	\delta_2 < \infty$. Here $c, \delta_1$ and $\delta_2$ are constants not depending on $n$. 
\end{Assumption}
\noindent
This assumption provides mild restrictions on the hyperparameters for the 
DAG-Wishart distribution. The assumption $2 < \alpha_i(\mathscr D) - 
\nu_i(\mathscr D)$ establishes prior propriety. The assumption $\alpha_i(\mathscr D) - 
\nu_i(\mathscr D) < c$ implies that the shape parameter $\alpha_i(\mathscr D)$ can only differ 
from $\nu_i(\mathscr D)$ (number of parents of $i$ in $\mathscr D$) by a constant which 
does not vary with $n$. Additionally, the eigenvalues of the 
scale matrix $U_n$ are assumed to be uniformly bounded in $n$.

For the rest of this paper, $p_n, {\Omega}_0^n, \Sigma_0^n, L_0^n, D_0^n, \mathscr{D}_0^n, \mathscr{D}^n,  d_n, q_n, \beta_n, \gamma_n, \tau_n, A_n$ will be denoted as $p, {\Omega}_0, \Sigma_0, L_0,$ $D_0, \mathscr{D}_0, \mathscr{D}, d, q, \beta, \gamma, \tau, A$ as needed for notational convenience and ease of exposition. We now state and prove the main joint variable and DAG selection consistency results.

\subsection{Posterior ratio consistency of $\gamma$ and $\mathscr D$}

In this section, we show that our method guarantees the posterior ratio consistency of $\gamma$ and $\mathscr D$.
Although \citet{Joint:Network:2016} consider a similar network-structured regression model, theoretical properties of Bayesian models such as  posterior ratio consistency and joint selection consistency have not been established yet up to our knowledge.
We first establish the posterior ratio consistency with respect to $\mathscr D$ under the true variable indicator $\gamma_0$.
Theorem \ref{thm1} says that the true DAG will be the mode of the posterior distribution with probability tending to $1$ as $n \rightarrow \infty$ under fixed $\gamma_0$.

\begin{theorem}\label{thm1}
	Under  Assumptions \ref{assumption:precision} \ref{assumption:subexponetial}, \ref{assumption:q} and \ref{assumption:hyper},
	$$\max_{\mathscr D \neq \mathscr D_0}\frac{\pi(\gamma_0, \mathscr D|   Y , X)}{\pi(\gamma_0, \mathscr D_0|   Y, X)}  \stackrel{\bar{P}}{\rightarrow} 0, \quad \mbox{as } n \rightarrow \infty.$$
\end{theorem}

\begin{remark}\label{remark:ordering:DAG}
	We would like to point out that the posterior ratio consistency for DAG is achieved under a given parent ordering, where that all the edges are directed from larger vertices to smaller vertices. For several applications in 
	genetics and environmental sciences, a location or time based ordering of variables 
	is naturally available. For temporal data, a natural ordering of variables is provided by 
	the time at which they are observed. In quantitative molecular applications,  the variables can be 
	genes or SNPs located on a chromosome, and their spatial location 
	provides a natural ordering. See \citet{HLPL:2006, Shojaie:Michailidis:2010, Yu:Bien:2016, KORR:2017} and references therein. 
\end{remark}

The next theorem establishes the posterior ratio consistency with respect to $\gamma$ under DAG $\mathscr D$. 
This notion of consistency implies that the true variable indicator $\gamma_0$ will be the mode of the posterior distribution with probability tending to $1$ as $n \rightarrow \infty$ under fixed $\mathscr D$. 
\begin{theorem} \label{thm2}
	Under Assumptions \ref{assumption:precision}-\ref{assumption:hyper}, the following holds: 
	$$\max_{(\gamma, {\mathscr{D}}) \ne (\gamma_0, {\mathscr{D}}_0)} \frac{\pi(\gamma, \mathscr D|   Y, X)}{\pi(\gamma_0, \mathscr D|   Y, X)}  \stackrel{\bar{P}}{\rightarrow} 0, \quad \mbox{as } n \rightarrow \infty.$$
\end{theorem}

\begin{remark} \label{remark:ordering:gamma}
	Based on a reviewer's comment, by carefully examining the proof of Theorem \ref{thm2}, we find out that even under a DAG with mis-specified ordering, the consistency result for $\gamma$ under fixed $\mathscr D$ will still hold. 
	We also investigate the performance of the proposed method under mis-specified ordering in Section \ref{sec:experiments}.
	The results suggest that our method recovers the true variable indicator $\gamma_0$ well even in the mis-specified case.   
\end{remark}


\noindent
From Theorem \ref{thm1}, Theorem \ref{thm2} and the fact that 
\begin{eqnarray*}
	\frac{\pi(\gamma, \mathscr D|   Y, X)}{\pi(   \gamma_0, \mathscr D_0|   Y, X)} &=& \frac{\pi(   \gamma_0, \mathscr D|   Y, X)}{\pi(   \gamma_0, \mathscr D_0|   Y, X)} \times \frac{\pi(   \gamma, \mathscr D|   Y, X)}{\pi(\gamma_0, \mathscr D|   Y, X)} ,
\end{eqnarray*}
we can obtain the joint posterior ratio consistency with respect to both $\gamma$ and $\mathscr D$.
It implies that the true variable indicator and DAG, $(\gamma_0, \mathscr D_0)$,  will be the mode of the posterior distribution with probability tending to $1$.
\begin{theorem} \label{thm3}
	Under Assumptions \ref{assumption:precision}-\ref{assumption:hyper}, the following holds: 
	$$\max_{(\gamma, {\mathscr{D}}) \ne (\gamma_0, {\mathscr{D}}_0)}\frac{\pi(\gamma, \mathscr D|   Y, X)}{\pi(\gamma_0, \mathscr D_0|   Y, X)}  \stackrel{\bar{P}}{\rightarrow} 0 \quad \mbox{as } n \rightarrow \infty,$$ 
	which implies that 
	$$\bar{P}((\hat{\gamma}, \hat{\mathscr D}) = (\gamma_0, \mathscr D_0)) \rightarrow 1, \mbox{ as } n \rightarrow \infty.$$ 
\end{theorem}

\subsection{Strong selection consistency of $\gamma$ and $\mathscr D$}

In this section, we establish the joint strong selection consistency with respect to both $\gamma$ and $\mathscr D$. 
Theorem \ref{thm4} shows that the posterior probability assigned to the true variable indicator $\gamma_0$ and the true underlying graph $\mathscr D_0$ grows to $1$ as $n \rightarrow \infty$.
We call this property the joint strong selection consistency.
Note that the result given in Theorem \ref{thm3} does not guarantee the joint strong selection consistency.

\begin{theorem}\label{thm4}
	Under Assumptions \ref{assumption:precision}-\ref{assumption:hyper}, if we further assume $\alpha_1 > 2$, the following holds: 
	$${\pi(\gamma_0, \mathscr D_0|   Y, X)}  \stackrel{\bar{P}}{\rightarrow} 1 \quad \mbox{as } n \rightarrow \infty.$$
\end{theorem}

We would like to point out that the condition on $\alpha_1$, which controls the rate of independent Bernoulli probability specified in Assumption \ref{assumption:q}, is only needed for strong selection consistency (Theorem \ref{thm4}). 
Similar restrictions on the hyperparameters have been considered for establishing consistency properties in the regression setup \citep{yang:2016, LLL:2018, CKG:nonlocal}. 
The model selection consistency for the posterior mode in Theorem \ref{thm3} does not require any restriction on $\alpha_1$.

All the aforementioned theorems are based on known $\sigma^2$, which tends to be not flexible enough, as in real applications, the underlying true variance often remains unavailable. Therefore, we introduce the following corollary for a fully Bayesian hierarchical approach, where an appropriate inverse-gamma prior is imposed on $\sigma^2$.
It turns out that even with the unknown $\sigma^2$, strong model selection consistency still holds under the same conditions given in Theorem \ref{thm4}.
\begin{corollary} \label{cor1}
	Suppose $\sigma^2$ is unknown and a proper inverse-gamma density with some positive constant parameters $(a_0, b_0)$ is placed on $\sigma^2$. Under Assumptions \ref{assumption:precision}-\ref{assumption:hyper}, and $\alpha_1 > 2$, the following holds: 
	$${\pi(\gamma_0, \mathscr D_0|   Y, X)}  \stackrel{\bar{P}}{\rightarrow} 1 \quad \mbox{as } n \rightarrow \infty.$$
\end{corollary}

\subsection{Behavior of the posterior probability when $b = 0$} \label{sec:b=0}
In this section, we aim to examine the behavior of the posterior probability for $(\gamma_0,\mathscr D_0)$ corresponding to two different scenarios when the MRF prior parameter $b > 0$ and $b = 0$ respectively. The goal is to show that under certain assumption on the connection between the sparsity patterns in $\gamma_0$ and $\mathscr D_0$, by borrowing the graph information through the MRF prior, the posterior probability assigned to $(\gamma_0,\mathscr D_0)$ will increase. In particular, we introduce the following condition with respect to the true sparsity patterns encoded in both the variable indicator and the graph.

\begin{condition} \label{condition:A}
	The true adjacency matrix $G_0$ and the true variable indicator $\gamma_0$ satisfy ${\gamma_0}_i = {\gamma_0}_j = 1$ whenever $(G_0)_{ij} = 1$ for $1 \le i, j \le p$.
\end{condition}
\noindent
Condition \ref{condition:A} essentially assumes that the connected variables through the underlying true DAG are active.
Under this condition, compared with modeling the variable and DAG separately, i.e. $b = 0$, incorporating network information into variable selection through the MRF prior with $b > 0$ will increase the posterior probability assigned to $(\gamma_0,\mathscr D_0)$ as illustrated in the following theorem. Proof for Theorem \ref{thm:b=0} will again be provided in the Supplementary material.

\begin{theorem} \label{thm:b=0}
	Let $\pi_1(\gamma_0,\mathscr D_0 \mid Y,X)$ be the posterior probability evaluated at $(\gamma_0,\mathscr D_0)$ under $b > 0$ and $\pi_2(\gamma_0,\mathscr D_0 \mid Y,X)$ be the posterior probability evaluated at $(\gamma_0,\mathscr D_0)$ under $b = 0$. 
	The following holds:
	$$\pi_1(\gamma_0,\mathscr D_0 \mid Y,X) > \pi_2(\gamma_0,\mathscr D_0 \mid Y,X).$$
\end{theorem}
\noindent
Theorem \ref{thm:b=0} implies that, under Condition \ref{condition:A}, our method achieves joint strong selection consistency without the condition on $b$ stated in Assumption \ref{assumption:hyper:markov}, which means the hyperparameter $b$ in the MRF prior does not need to go to zero.

\setcounter{equation}{0} 
\section{Numerical Studies}\label{sec:experiments}

\subsection{Posterior inference}\label{subsec:postinf}
For given positive real values $a_0$ and $b_0>0$, let $IG(a_0, b_0)$ be the inverse-gamma distribution with the shape parameter $a_0$ and scale parameter $b_0$.
Then, similar to \eqref{conditional}, the joint posterior distribution of $\gamma$ and $\mathscr D$ based on \eqref{model:spec1}--\eqref{model:spec6} and $\sigma^2 \sim IG(a_0, b_0)$ is 
\bea
&&  \pi(\gamma , \mathscr D \mid Y, X )  \\
&\propto& \pi(\gamma|\mathscr D)\pi(\mathscr D) \frac{z_{\mathscr{D}}(U+X^TX,n+\alpha(\mathscr D))}{z_{\mathscr{D}}(U, \alpha(\mathscr D))} \\
&& \times \,\, \det\left(I_{|\gamma|} + \tau^2X_\gamma^T X_\gamma \right)^{-\frac 1 2 } \Big\{  b_0 + \frac{1}{2} Y^T \left( I_n + \tau^2 X_\gamma X_\gamma^T \right)Y \Big\}^{- \frac{n+2a_0}{2}}  .
\eea
We suggest using a Metropolis-Hastings within Gibbs sampling for posterior inference: 
\begin{enumerate}
	\item Set the initial values $\gamma^{(1)}$ and $\mathscr D^{(1)}$. \vspace{-.3cm}
	\item For each $s=2,\ldots, S$,  \vspace{-.3cm}
	\begin{enumerate}
		\item sample $\gamma^{new} \sim q_\gamma(\cdot \mid \gamma^{(s-1)})$;
		\item set $\gamma^{(s)}= \gamma^{new}$ with the probability 
		\bea
		p_{acc,\gamma} &=& \min \left\{ 1,  \frac{\pi(\gamma^{new} \mid \mathscr D^{(s-1)}, Y,X)}{\pi(\gamma^{(s-1)} \mid \mathscr D^{(s-1)}, Y,X)}  \frac{q_\gamma(\gamma^{(s-1)} \mid \gamma^{new}) }{q_\gamma(\gamma^{new} \mid \gamma^{(s-1)})}  \right\},
		\eea
		otherwise set $\gamma^{(s)}= \gamma^{(s-1)}$;
		\item sample $\mathscr D^{new} \sim q_{\mathscr D}(\cdot \mid \mathscr D^{(s-1)})$;
		\item set $\mathscr D^{(s)}=\mathscr D^{new}$ with the probability 
		\bea
		p_{acc,\mathscr D} &=& \min \left\{ 1,  \frac{\pi(\mathscr D^{new} \mid \gamma^{(s)}, Y,X)}{\pi(\mathscr D^{(s-1)} \mid \gamma^{(s)}, Y,X)}  \frac{q_{\mathscr D}(\mathscr D^{(s-1)} \mid \mathscr D^{new}) }{q_{\mathscr D}(\mathscr D^{new} \mid \mathscr D^{(s-1)})}  \right\},
		\eea
		otherwise set $\mathscr D^{(s)}= \mathscr D^{(s-1)}$.
	\end{enumerate}
\end{enumerate}
The inference for the DAG $\mathscr D$, the steps 2-(c) and 2-(d) in the above algorithm, can be parallelized for each column.
For more details, we refer to \cite{CKG:2017} and \cite{LLL:2018}.
We used the proposal kernel $q_\gamma(\cdot \mid \gamma')$ which gives a new set $\gamma^{new}$ by changing a randomly chosen nonzero component in $\gamma'$ to $0$  with probability $0.5$ or by changing a randomly chosen zero component to $1$ randomly with probability $0.5$.
The same proposal kernels were used for each column of $\mathscr D$.

\subsection{Simulation Studies}\label{subsec:simulation}

In this section, we demonstrate the performance of the proposed method in various settings.
We closely follow but slightly modify the simulation settings in \cite{Joint:Network:2016}.

Suppose that we have $X_i = (X_{i1}, \ldots, X_{ip})^T \overset{i.i.d.}{\sim} N_p(0, \Sigma_0), i=1,\ldots, n$, where $\Sigma_0^{-1} = L_0 (D_0)^{-1} L_0^T$, $n=100$ and $p=240$.
If we consider $p$ as the number of genes, we have $240$ genes in this case.
Among $240$ genes, we assume that there are $40$ transcription factors (TFs) and each TF regulates $5$ genes.
Let $TF_j$ be the index for the $j$th TF and $(TF_1,  TF_2,\ldots, TF_{40}) = (6, 12, \ldots, 240)$.
It corresponds to the DAG $\mathscr D_0$, the support of $L_0$, such that $pa_{TF_j - k}(\mathscr D_0) = \{  TF_j \}$ for $j=1,\ldots, 40$ and $k=1,\ldots, 5$.
Suppose that the TFs independently follow the normal distribution, that is, $X_{TF_j} \overset{ind}{\sim}  N(0,d_{TF_j})$, where $d_{TF_j} \overset{i.i.d.}{\sim} Unif(3,5)$, for $j=1,\ldots, 40$.
We further assume that, given $X_{TF_j}$, the conditional distribution of the gene $X_j$ that $TF_{j'}$  regulates is $N(  X_{TF_{j'}} , d_j)$, where $d_j \overset{i.i.d.}{\sim} Unif(3,5)$ for $j=1,\ldots, 240$.
It corresponds to the true modified Cholesky parameter $(L_0, D_0)$ such that $(L_0)_{TF_j, TF_j - k} =1$ and $D_0 = diag(d_j)$ for $j=1,\ldots, 40$ and $k=1,\ldots, 5$.
We simulate the data from
\bea
Y &=& X\beta_0 +\epsilon,
\eea
where $X = (X_1,\ldots, X_n)^T$ and $\epsilon \sim N_n(0,  \sigma_\epsilon^2 I_n )$ and $\sigma_\epsilon^2 = \|\beta_0\|_2^2/4$.
We investigate four settings for the true coefficient vector $\beta_0$ as described in \citet{li2008network} and \cite{Joint:Network:2016}.
In the first setting, it is assumed that $\beta_{0,TF_{1:4}} = (5,-5, 3,-3)^T$,  $\beta_{0,TF_{j}-k} = \beta_{0,TF_{j}} /\sqrt{10}$ for $j=1,2,3,4$ and $k=1,\ldots,5$, and $\beta_{0,j} = 0$ for $j=25,\ldots, 240$.
This setting implies that the genes in the same cluster have the same signs for the coefficients.
In the second setting, the true coefficient $\beta_0$ is the same as the first setting except that the signs are reversed for the two genes that $TF_j$ regulates, i.e., $\beta_{0,TF_{j}-k} = -\beta_{0,TF_{j}} /\sqrt{10}$ for $j=1,2,3,4$ and $k=1,2$.
This setting implies that the genes in the same cluster might have different signs for the coefficients.
The third and fourth settings the same as the first and second settings expect considering $10$ instead of $\sqrt{10}$.
Thus, they consider smaller signals.
We call this simulation setting Scenario 1.

We also investigate a different simulation scenario, say Scenario 2, where the signals in $\beta_0$ are small. 
In this case, there are $p=150$ genes, $30$ TFs and $4$ regularized genes for each TF.
The precision matrix $\Sigma_0^{-1} = L_0 (D_0)^{-1} L_0^T$ is generated by $d_j \overset{i.i.d.}{\sim} Unif(2,5)$ and $(L_0)_{TF_j, TF_j - k} \overset{i.i.d.}{\sim} Unif(0.3, 0.7)$.
The variance of $\epsilon$ is chosen as $\sigma_\epsilon^2 = \|\beta_0\|_2^2$.
We consider four settings for the true coefficient vector $\beta_0$.
In the first and third settings, $\beta_0$ is generated by $\beta_{0,j} \overset{i.i.d.}{\sim} Unif(0.5,1)$ and $\beta_{0,j} \overset{i.i.d.}{\sim} Unif(0.2,1)$ for $j=1,\ldots, 20$, respectively, and $\beta_{0,j}=0$ for $j=21,\ldots, 150$.
In the second and fourth settings, we only change the signs of nonzero entries of $\beta_0$ randomly.
We call this simulation setting Scenario 2.

Lastly, we consider a setting where the network structure of the covariate $X$ is a undirected graph. 
We generate the covariates $\tilde{X}_i \overset{i.i.d.}{\sim} N_p(0 , \sg_0 ), i=1,\ldots,n$, where $n=100, p=150$, $\sg_0 = \tilde{\sg}_0 + \{0.01 - eig_1(\tilde{\sg}_0)\} I_p$ and  
\bea
(\tilde{\sg}_0)_{ij} = 
\begin{cases}
	2 \max \big(  1 - \frac{|i-j|}{10} , 0 \big), & \text{ if } |i-j| \le 5\\
	0, &  \text {otherwise}.
\end{cases}
\eea
Note that $\sg_0$ is positive definite.
Furthermore, to consider the mis-specified ordering case, we randomly shuffle columns of $\tilde{X} = (\tilde{X}_1, \ldots, \tilde{X}_n)^T$ to construct $X$.
We simulate the data from $Y= X\beta_0 + \epsilon$, where $\epsilon \sim N_n (0, \sigma_\epsilon^2 I_n)$ and $\sigma_\epsilon^2 = \|\beta_0\|_2^2/4$.
Two settings for the true coefficient vector $\beta_0$ are considered.
In the first setting, $\beta_0$ is generated by $\beta_{0,j} \overset{i.i.d.}{\sim} Unif(0.5, 1)$ for $j=1,\ldots, 10$ and $\beta_{0,j}=0$ for $j=11,\ldots,150$.
In the second setting, we only change the signs of nonzero entries of $\beta_0$ randomly.
We call this simulation setting Scenario 3, and the simulation results for this setting are reported at Table \ref{table:comp3}.

We compare the performance of our joint selection method with other existing variable selection methods: Lasso \citep{tibshirani1996regression}, elastic net \citep{zou2005regularization} and the Bayesian joint selection method proposed by \cite{Joint:Network:2016}.
The tuning parameters in Lasso and elastic net were chosen by 10-fold cross-validation.
For Bayesian methods, as discussed by \cite{Joint:Network:2016}, we suggest using the hyperparameters $a=2.75$ and $b=0.5$ for the MRF prior as default.
Furthermore, to show the benefits of joint modeling, we also tried the setting with $b=0$ which corresponds to the Bayesian method modeling the variable
and DAG separately.
The other hyperparameters were set at $a_0=0.1, b_0= 0.01,\tau^2= 1, q=0.005, U = I_p$ and  $\alpha_i(\mathcal{D}) = \nu_i(\mathcal{D})+10$ for all $i=1,\ldots, p$.
The initial state for $\gamma$ was set at $p$-dimensional zero vector, i.e., the empty model, while the initial state for $\mathscr D$ was chosen by the CSCS method \citep{KORR:2017}.
For posterior inference, $5,000$ posterior samples were drawn after a burn-in period of $5,000$.
The indices having posterior inclusion probability larger than $0.5$ were included in the final model.
The resulting model is called the median probability model, and when there is a model with posterior probability larger than $1/2$, it coincides with the posterior mode \cite{barbieri2004optimal}.
Since we have proved the joint strong selection consistency (Theorem \ref{thm4}), the two models are asymptotically equivalent in our setting.
Thus, although other approaches (for example, see \cite{scott2008feature}) can be adapted to give a reasonable estimate of the posterior mode, we use the median probability model as a convenient but asymptotically equivalent alternative.

To evaluate the performance of variable selection, the sensitivity, specificity, area under the curve (AUC), Matthews correlation coefficient (MCC), the number of errors (\#Error) and mean-squared prediction error (MSPE) are reported at Tables \ref{table:comp1}, \ref{table:comp2} and \ref{table:comp3}.
The criteria are defined as
\bea
\text{Sensitivitiy}  &=&     \frac{TP}{TP+FN} ,   \\
\text{Specificity}  &=&    \frac{TN}{TN+FP}   ,  \\
\text{MCC}  &=&     \frac{TP \times TN - FP\times FN}{\sqrt{(TP+FP)(TP+FN)(TN+FP)(TN+FN)}}  ,	  \\
\#\text{Error} &=& FP + FN , \\
\text{MSPE}  &=&    \frac{1}{n_{\rm test}} \sum_{i=1}^{n_{\rm test}}     \big( \hat{Y}_i - Y_{{\rm test},i}  \big)^2  ,
\eea
where TP, TN, FP and FN are true positive, true negative, false positive and false negative, respectively.
The AUC is calculated based on the true positive rate (Sensitivity) and the false positive rate ($1-$Specificity) for Bayesian methods with varying thresholds.
To draw the AUC, for each threshold, the indices having posterior inclusion probability larger than a given threshold were included in the final model.
The AUCs for the regularization methods are omitted.
We denote $\hat{Y}_i = X_i^T \hat{\beta}$, where $\hat{\beta}$ is the estimated coefficient based on each method.
For Bayesian methods, the usual least square estimates based on the selected support were used as $\hat{\beta}$.
We generated  test samples $Y_{{\rm test}, 1},\ldots, Y_{{\rm test}, n_{\rm test}}$ with $n_{\rm test}=100$ to calculate the MSPE.

\begin{table}[!tb]
	\centering\tiny
	\caption{
		The summary statistics for Scenario 1 are represented for each setting.
		Different setting means different choice of the true coefficient $\beta_0$.
		Sens and Spec are sensitivity and specificity, respectively.
		Joint.CL: the Bayesian joint selection method proposed in this paper. 
		Joint.P: the Bayesian joint selection method suggested by \cite{Joint:Network:2016}.
		Elastic: elastic net.
	}\vspace{.15cm}
	\begin{tabular}{c c c c c c c | c c c c c c}
		\hline 
		& \multicolumn{6}{c}{ Setting 1 } & \multicolumn{6}{c}{ Setting 2 } \\ 
		& Sens & Spec & AUC & MCC & \#Error & MSPE & Sens & Spec & AUC & MCC & \#Error & MSPE \\ \hline
		Joint.CL $(b=\frac{1}{2})$ & 0.8750 & 0.9861 & 0.9937  &   0.8611 & 6  & 69.1445  & 0.8750 & 0.9954 & 0.9894   & 0.9049  &  4  & 56.4885   \\ 
		Joint.CL $(b=0)$ & 0.7500 & 0.9815 & 0.9601 & 0.7605  &  10  & 96.9889  & 0.3333 & 1.0000 & 0.9058   & 0.5571    &  16  &  142.2708   \\ 
		Joint.P  & 0.8750 & 0.9861 & 0.9838   &   0.8611 & 6  & 71.0443 &   0.7500 & 0.9954 & 0.9958   & 0.8282  &  7  & 73.7870  \\ 
		Lasso  & 1.0000 & 0.8056 & $\cdot$   &  0.5412   &  42  & 45.5522  &  0.7083 & 0.8519 & $\cdot$      & 0.4170  &   39   &  106.0526  \\ 
		Elastic & 1.0000 & 0.9352 & $\cdot$    &  0.7685  &   14  &  41.8631 &  0.8750 & 0.8426 & $\cdot$      & 0.5122   &  37   &  92.6665  \\ \hline \hline 
		& \multicolumn{6}{c}{ Setting 3 } & \multicolumn{6}{c}{ Setting 4 } \\ 
		& Sens & Spec & AUC & MCC & \#Error & MSPE & Sens & Spec & AUC & MCC & \#Error & MSPE \\ \hline
		Joint.CL $(b=\frac{1}{2})$ & 0.2083 & 0.9907 & 0.8493   & 0.3549  & 21   & 42.5213 & 0.3750 & 1.0000 & 0.7373   & 0.5922    &  15  &  30.3394  \\ 
		Joint.CL $(b=0)$ & 0.1667 & 0.9907  & 0.7117  & 0.3025  & 22   & 42.7116 & 0.3333 & 0.9954 & 0.7619   & 0.5191    &  17  &  35.0479  \\ 
		Joint.P  & 0.2500 & 0.9907 & 0.8559   &   0.4023 & 20  & 40.3569 &   0.2917 &  0.9954  & 0.8954  &  0.4797  &  18  & 35.7181  \\ 
		Lasso  & 1.0000 & 0.8241 & $\cdot$      &  0.5648   &  38  & 32.1919  &  0.6667 & 0.8102  & $\cdot$     & 0.3362  &   49   &  40.7437  \\ 
		Elastic & 1.0000 & 0.9444  & $\cdot$      &  0.7935  &  12   &  29.3908 &  0.6250 & 0.8935  & $\cdot$     & 0.4261   &  32   &  34.9673  \\ \hline  
	\end{tabular}\label{table:comp1}
\end{table}

\begin{table}[!tb]
	\centering\tiny
	\caption{
		The summary statistics for Scenario 2 are represented for each setting.
		Different setting means different choice of the true coefficient $\beta_0$.
	}\vspace{.15cm}
	\begin{tabular}{c c c c c c c | c c c c c c}
		\hline 
		& \multicolumn{6}{c}{ Setting 1 } & \multicolumn{6}{c}{ Setting 2 } \\ 
		& Sens & Spec & AUC & MCC & \#Error & MSPE & Sens & Spec & AUC & MCC & \#Error & MSPE \\ \hline
		Joint.CL $(b=\frac{1}{2})$ & 0.7500 & 0.9923 & 0.9362  & 0.8174  &  6  & 20.9925  & 0.6000 & 1.0000 & 0.8933    & 0.7518   &  8  & 15.8789   \\ 
		Joint.CL $(b=0)$ & 0.6000 & 1.0000 & 0.9200    & 0.7518  & 8  & 29.6691 & 0.6500 & 0.9923 & 0.8790    & 0.7506    &  8  &  23.3007  \\ 
		Joint.P  & 0.6500 & 1.0000 & 0.9842   & 0.7854 & 7  & 15.4705 &   0.5000 & 1.0000 & 0.9081    & 0.6814   &  10  &  19.2450  \\ 
		Lasso  & 1.0000 & 0.8308 & $\cdot$    &  0.6290   &  22  & 14.8092 &  0.9000 & 0.7692 & $\cdot$    & 0.4877  &   32   &  13.4260  \\ 
		Elastic & 0.9500 & 0.9077  & $\cdot$    &  0.7201  &   13  &  18.9942 &  0.8000 & 0.8615 & $\cdot$    & 0.5371   &  22   &  14.5779 \\ \hline \hline 
		& \multicolumn{6}{c}{ Setting 3 } & \multicolumn{6}{c}{ Setting 4 } \\ 
		& Sens & Spec & AUC & MCC & \#Error & MSPE & Sens & Spec & AUC & MCC & \#Error & MSPE \\ \hline
		Joint.CL $(b=\frac{1}{2})$ & 0.7500 & 1.0000 & 0.9537    & 0.8498  & 5   & 6.6246 & 0.6500 & 1.0000   & 0.8398  & 0.7854    &  7  &  7.4111  \\ 
		Joint.CL $(b=0)$ & 0.4000 & 1.0000  & 0.9631   & 0.6051  & 12   & 20.3681 & 0.3000 & 1.0000  & 0.7962   & 0.5204    &  14  &  12.7521  \\ 
		Joint.P  & 0.6500 & 1.0000 & 0.9811    &   0.7854 & 7  &  11.6528 &   0.4500 & 1.0000 & 0.9057    & 0.6441    &  11  &  9.2049  \\ 
		Lasso  & 0.9500 & 0.8154 & $\cdot$    &  0.5754  &  25  & 8.2451  &  0.8500 & 0.7462 & $\cdot$    & 0.4299  &   36   &  7.7223  \\ 
		Elastic & 0.9500 & 0.8923 & $\cdot$    &  0.6912  &  15   &  10.7742 &  0.7000 & 0.8846 & $\cdot$   & 0.5032   &  21   &  7.8241 \\ \hline  
	\end{tabular}\label{table:comp2}
\end{table}

Based on Tables \ref{table:comp1} and \ref{table:comp2}, we notice that Bayesian joint selection methods tend to have better specificity and MCC, while the regularization methods (Lasso and elastic net) have better sensitivity.
As discussed by \cite{Joint:Network:2016}, this seems natural because the regularization methods based on cross-validation tend to include many redundant variables.
It leads to relatively larger number of  errors for the regularization methods compared with those for the Bayesian joint selection methods.
We also found that the joint Bayesian  selection method proposed in this paper (Joint.CL $(b=1/2)$) works better than that proposed by \cite{Joint:Network:2016} (Joint.P) in terms of performance measures in Tables \ref{table:comp1} and \ref{table:comp2} except the AUC. 
In fact, the two Bayesian joint selection methods are quite similar to each other except the graph structure they consider.
In these simulation scenarios, the DAG structure seems more appropriate because clearly there are parents (TFs genes) and children (regularized genes for each TF).
Thus, our method would be preferable in this case.
Furthermore, based on asymptotic results, one can expect that our method will give accurate inference results as we have more observations, while asymptotic properties of the Bayesian method proposed by \cite{Joint:Network:2016}  are still in question. 
Lastly, the results show that our joint modeling (Joint.CL $(b=1/2)$) significantly improves the performance of variable selection compared with modeling the variable and DAG separately (Joint.CL $(b=0)$).
These suggest that joint modeling approach actually improves the performance of variable selection by borrowing information from the DAG structure.

\begin{table}[!tb]
	\centering\tiny
	\caption{
		The summary statistics for Scenario 3 are represented for each setting.
		Different setting means different choice of the true coefficient $\beta_0$.
	}\vspace{.15cm}
	\begin{tabular}{c c c c c c c | c c c c c c}
		\hline 
		& \multicolumn{6}{c}{ Setting 1 } & \multicolumn{6}{c}{ Setting 2 } \\ 
		& Sens & Spec & AUC & MCC & \#Error & MSPE & Sens & Spec & AUC & MCC & \#Error & MSPE \\ \hline
		Joint.CL $(b=\frac{1}{2})$ & 1.0000 & 0.9357 & 0.9964  & 0.7018  &  9  & 1.6875  & 1.0000 & 0.9500 & 0.9821    & 0.7475   &  7  & 1.6469   \\ 
		Joint.CL $(b=0)$ & 1.0000 & 0.9429 & 0.9786    & 0.7237  & 8  & 1.7331 & 1.0000 & 0.9429 & 0.9786   & 0.7237    &  8  &  1.7331  \\ 
		Joint.P  & 1.0000 & 0.9500 & 0.9857   & 0.7475 & 7  & 1.6838 &   1.0000 & 0.9500 & 0.9821    & 0.7475   &  7  &  1.6838  \\ 
		Lasso  & 1.0000 & 0.3571 & $\cdot$    &  0.1890   &  90  & 1.8121 &  1.0000 & 0.3571 & $\cdot$    & 0.1890  &   90   &  1.8121  \\ 
		Elastic & 1.0000 & 0.6714  & $\cdot$    &  0.3463  &   46  &  1.6770 &  1.0000 & 0.6286 & $\cdot$    & 0.3184   &  52   &  1.6969 \\ \hline 
	\end{tabular}\label{table:comp3}
\end{table}
Table \ref{table:comp3} shows the results for Scenario 3, where the true network structure for $X$ is a undirected graph and the ordering is mis-specified.
Even in this case, our joint modeling method provides comparable performance to that of \cite{Joint:Network:2016}, which is designed for undirected graphs.
Similar to Scenarios 1 and 2, regularization methods do not work well compared with Bayesian methods in our settings.

\setcounter{equation}{0} 
\section{Discussion} \label{sec:discussion}
In this paper, we work in a regression setting, where the predictors are both relevant to a response variable of interest and  functionally related to one another via a Gaussian DAG model. In particular, we consider a hierarchical multivariate regression model with DAG-Wishart priors on the covariance matrix for the predictors, spike and slab priors on regression coefficients, independent Bernoulli priors for each edge in the DAG, and a MRF prior linking the variable indicators to the graph structure. Under high-dimensional settings and standard regularity assumptions, when the underlying variance $\sigma^2$ is available, we establish both posterior ratio consistency and strong selection consistency for estimating the variable and the graph for the covariates jointly. When the underlying response variance is unknown and an appropriate inverse gamma prior is placed on $\sigma^2$, we also establish the joint selection consistency under the same regularity conditions. Finally, through simulation studies, we demonstrate that the model studied in this paper can outperform existing state-of-the-art methods in selecting network-structured predictors including both penalized likelihood and Bayesian approaches in several settings. For future studies, we intend to explore other types of priors over the graph space and on the regression coefficients to see if the consistency and better simulation performance can both be achieved under weakened assumptions.

\vskip 14pt
\noindent {\large\bf Supplementary Materials}\\
Supplementary material includes the proofs for main results and other auxiliary results.
\par
\vskip 14pt
\noindent {\large\bf Acknowledgements}\\
We thank Dr. Christine Peterson for sharing the code to implement joint Bayesian variable and graph selection method in \cite{Joint:Network:2016}.
We would like to thank two referees for their valuable comments which have led to improvements of an earlier version of the paper. 
This research was supported by the Simons Foundation's collaboration grant (No.635213), the National Research Foundation of Korea (NRF) grant funded by the Korea government (MSIT) (No.2019R1F1A1059483) and INHA UNIVERSITY Research Grant.
\par

\markboth{\hfill{\footnotesize\rm XUAN CAO AND KYOUNGJAE LEE} \hfill}
{\hfill {\footnotesize\rm BAYESIAN JOINT SELECTION} \hfill}

\bibhang=1.7pc
\bibsep=2pt
\fontsize{9}{14pt plus.8pt minus .6pt}\selectfont
\renewcommand\bibname{\large \bf References}
\expandafter\ifx\csname
natexlab\endcsname\relax\def\natexlab#1{#1}\fi
\expandafter\ifx\csname url\endcsname\relax
  \def\url#1{\texttt{#1}}\fi
\expandafter\ifx\csname urlprefix\endcsname\relax\def\urlprefix{URL}\fi

\bibliographystyle{chicago}      
\bibliography{references}   

\vskip .65cm
\noindent
University of Cincinnati
\vskip 2pt
\noindent
E-mail: (caox4@ucmail.uc.edu)
\vskip 2pt

\noindent
Inha University
\vskip 2pt
\noindent
E-mail: (leekjstat@gmail.com)

\end{document}